\documentclass{article}

\usepackage{amssymb}

\newcommand{\NN}{\mathbb{N}}
\newcommand{\PP}{\mathbb{P}}
\newcommand{\RR}{\mathbb{R}}
\newcommand{\F}{\mathbb{F}}
\newcommand{\SaS}{\mathbb{S}}

\newcommand{\T}{\mathbb{T}}

\newcommand{\D}{{\cal D}}

\newcommand{\K}{{\cal K}}
\newcommand{\Aa}{{\cal A}}

\newcommand{\Le}{{\cal L}}
\newcommand{\hH}{{\cal H}}
\newcommand{\I}{{\cal I}}
\newcommand{\eE}{{\cal E}}

\newcommand{\Y}{{\cal Y}}
\newcommand{\wY}{{\widehat{Y}}}
\newcommand{\tY}{{\underline{Y}}}

\newcommand{\R}{{\Theta}}
\newcommand{\aS}{{\cal S}}

\newcommand{\ve}{{\vec e}}

\newcommand{\htau}{\underline{\tau}}
\newcommand{\wtau}{\widehat{\tau}}
\newcommand{\hGamma}{\widehat{\Gamma}}

\newcommand{\Int}{{\rm Int}\,}
\newcommand{\qperp}{{\underline{\perp}}}

\newcommand{\hLambda}{\widehat{\Lambda}}

\newtheorem{theorem}{Theorem}[section]
\newtheorem{remark}{Remark}[section]
\newtheorem{proposition}{Proposition}[section]

\begin{document}

\begin{center}
{\Large STIT Process and Trees}

\vspace{1cm}

{Servet Mart{\'i}nez}\\
{Departamento Ingenier{\'\i}a Matem\'atica and Centro
Modelamiento Matem\'atico,\\  Universidad de Chile,\\
UMI 2807 CNRS, Casilla 170-3, Correo 3, Santiago, Chile.\\
Email: smartine@dim.uchile.cl} 

\vspace{0.5cm}

{Werner Nagel}\\
{Friedrich-Schiller-Universit\"at Jena,\\
Institut f\"ur Stochastik,\\
Ernst-Abbe-Platz 2,
D-07743 Jena, Germany.\\
Email: werner.nagel@uni-jena.de} 

\end{center}

\begin{abstract}
We study several constructions of
the STIT tessellation process in a window of $\RR^\ell$
and supply an exact formula for its transition probability.
\end{abstract}

\noindent {\em Keywords:}
Stochastic geometry, Random process of tessellations, 
STIT tessellation, Binary tree


\vspace{1cm}

\noindent {\em AMS subject classification:} 60D05; 60J25; 60J75; 05C05

\section{Introduction}
\label{section0}

We will describe the STIT tessellation process $Y\wedge W$ on a 
window $W$ in $\RR^\ell$ as it was defined for the first time in \cite{nagel/weiss05}.

\medskip

A STIT process is a particular cell division process, and within a 
bounded window it is a pure jump Markov process. Hence, 
it can be considered from two 
aspects. One aspect is
that each cell has a random lifetime, and at the end of its lifetime 
the cell is divided and two new cells are born. 
Since the lifetimes of the cells run simultaneously, 
this approach can appropriately be 
described using binary rooted trees where the nodes 
represent the cells of the tessellation. 
The other aspect is that the STIT process in a bounded 
window has a random holding time in a state and when this time 
elapsed it jumps into another state. This jump is performed
by first a random selection of a cell that has to be divided and 
then dividing this cell. In the present paper we consider both 
these aspects in detail and we relate them to each other.

\medskip

The main result is a new construction of $Y\wedge W$ 
using a sequence of random hyperplanes generated
with a random measure. This is done in Section \ref{section8}. 
This method allows to gain in the efficiency of simulation of 
STIT, because all the hyperplanes are used in the construction
in contrast to other constructions using a rejection method where one must attend
that a random hyperplane cuts a prescribed cell.

\medskip

Our construction uses a sequence of hyperplanes with a 
random distribution which depends on the current tessellation 
(and hence it has to be adapted after each cell division). 
This is different from the construction done in 
Section $4$ in \cite{m/n/w08a} that uses a Poisson point process of hyperplanes 
with the fixed intensity, but it requires to be corrected 
because the process underestimates the rate of apparition of hyperplanes
in STIT. On the other hand, in \cite{mecke10a}  
the mean length of segments is computed for mixtures of tessellations 
in the case $\ell=2$, 
and a Poisson tessellation process is constructed with 
an appropriate intensity measure. This construction differs from the 
one we make in the last paragraph of Section \ref{section8}, 
since our construction has a random intensity measure and holds 
for any dimension.

\medskip

We use rooted binary (dyadic) trees for the description of the STIT process 
in a window. This is natural since 
STIT is a cell division process. In Section \ref{section4} 
we define this class of trees, and define finite trees and its leaves as
appropriate graph objects in our study. 

\medskip

The construction of STIT is formally done in Section \ref{section5},
and the rooted binary tree helps to write a simple algorithm. 
The root of the tree represents the window and each node
stands for a cell which appears in the cell division
procedure. When a cell is
divided, the corresponding node in the tree has two children, representing
the two new cells and denoted by '+' or '-', each symbol indicating the 
half-space of the dividing hyperplane. The lifetimes of cells are independent and
exponentially distributed. Even if Proposition \ref{lawstit} and
Proposition \ref{propglobal} are not new,
we recall them to provide explicit conditions
characterizing the STIT process and to identify the stability-under-iteration 
property as one close to what is called branching property of a fragmentation 
chain in \cite{ber}.

\medskip

In Section \ref{section6} we use the tree representation of STIT to give 
a formula for the marginal distribution of
$(Y\wedge W)_t$ at a fixed time $t$ by considering all the possible
paths on a binary  rooted tree.

\medskip

In Section \ref{section7} we revisit the construction of STIT in a window,
 and
the order of choosing the random objects is kept:
first one selects the cell that will be broken, and conditioned to it 
one chooses the hyperplane that cuts it. We summarize the results
of this construction in Proposition \ref{primeas}, that
only serves to order the elements but with no novel elements. 
We emphasize that this construction is not optimal
for simulations because the time of retrieving such an hyperplane can be 
highly-time consuming since it is based on a rejection method. 

\medskip

In Section \ref{section8} we modify the above construction with 
a different order of choosing the random objects: first one generates a random
hyperplane with a distribution depending on the whole tessellation in the time
of up-dating.
Then the cell to be divided is chosen with equal probability among all the 
cells being intersected (this is formula (\ref{relcr2})). 
This is done in detail in Theorem \ref{igual1},
which to our view gives a novel approach to construct STIT. 
We point out that a relation similar to (\ref{relcr2}) appears in \cite{st} 
page $9$, but with a different hyperplane measure.

\medskip

In Section \ref{section1} we give some useful facts in probability, mainly
on conditional independence and Lebesgue probability spaces. The basic notions 
and notations for tessellations are supplied in Section \ref{section2} and in
Section \ref{section3} we summarize the main elements of the random law 
of random hyperplanes 
and supply its main properties which ensure that
it serves to define a STIT tessellation processes.

\section{Preliminaries on Probability}
\label{section1}
Let $(\Omega, {\cal B}(\Omega), \PP)$ be a probability space
which is the basis for the construction of all the random objects
we will use.

\medskip

To describe  relations among the random objects it is useful
to introduce some notions and notation. Let $(D,{\cal D})$
be a measurable space. $g:\Omega\to D$ is a random variable if 
$g^{-1}(\cal D)\subseteq {\cal B}$
and we put $\sigma(g)=g^{-1}(\cal D)$. If $h$ is another random variable
we put $h\in \sigma(g)$ if 
$\sigma(h)\subseteq \sigma(g)$. Here $g$ and $h$ can also
be countable sequences of random variables.

\medskip

If $g$ and $h$ are two random variables,
we write $g\, \qperp \,  h$ when
$\sigma(g)$ and $\sigma(h)$ are ($\PP-$)\hfill\break independent.
If $g$, $h$ are two random variables 
and $\Aa\subseteq {\cal B}$ is a sub $\sigma-$field,
we express by $(g\, \qperp \,  h) \, | \, \Aa$  that $g$
and $h$ are conditionally independent given $\Aa$, that is
$\PP(D'\cap D''\, | \, \Aa)=\PP(D'\, | \, \Aa)\PP(D''\, | \, \Aa)$
(a.s.) for all $D'\in \sigma(g)$, $D''\in \sigma(h)$.
Also if $z$ is a random variable
we put  $(g\, \qperp \,  h) \, | \, z$
for $(g\, \qperp \,  h) \, | \, \sigma(z)$.

\medskip

Let $g:\Omega\to D$ be a random variable and $F$ be a probability measure 
on $(D,\D)$, by $g\sim F$ we mean that $g$ is distributed as $F$.
If $z$ is a random variable and $F(z)$ is a random distribution
depending on $z$, we write
$g \,| \, z\sim F(z)$ to express that the conditional distribution of
$g$ given $\sigma(z)$ is $F(z)$.

\medskip

Let $(D,{\cal B}(D),P)$ be a probability space such that: $D$ is a complete 
separable metric space, ${\cal B}(D)$ is its Borel $\sigma-$field completed with 
respect to the probability measure $P$. Also assume that $P$ is non-atomic.
Then, $(D,{\cal B}(D),P)$ is a Lebesgue probability space, see \cite{delar}.
This means that $(D,{\cal B}(D),P)$ and
$([0,1],{\cal B}[0,1], \lambda)$ are isomorphic, where $\lambda$ is the Lebesgue 
measure. That is, there exists an isomorphism $v:[0,1]\to D$, which is
a bimeasurable function such that $\lambda(v^{-1}(A))=P(A)$
for all $A\in {\cal B}(D)$. 

\medskip

Let $(D,{\cal B}(D),P)$ be a Lebesgue probability space.
Let $Q$ be a probability measure equivalent to $P$ with Radon-Nikodym
derivate $f=dQ/dP>0$ $\, P-$a.s. Then
also $(D,{\cal B}(D),Q)$ is isomorphic to $([0,1],{\cal B}[0,1], \lambda)$.
So $(D,{\cal B}(D),P)$ and $(D,{\cal B}(D),Q)$ are isomorphic:
there exists a bimeasurable function $\Xi:D\to D$ that satisfies
$P(\Xi^{-1}(A))=Q(A)$ for all $A\in {\cal B}(D)$. 
We will mainly consider Lebesgue probability spaces.

\section{Tessellations}
\label{section2}


\medskip

For a set $B\subseteq \RR^\ell$ we denote respectively by $\partial B$ 
and $\Int B$ the {\it boundary} and the {\it interior} of $B$.
A {\it polytope} $K$ is the convex hull of a finite point set. 
A polytope with nonempty interior will be called {\it a cell} or 
{\it a window}. This distinction will depend on the context, 
usually, we reserve the name cell when the polytope 
belongs to a tessellation.

\medskip

A tessellation $T$ in $\RR^\ell$ is a locally finite class of cells with
disjoint interiors and covering the Euclidean space. 
The locally finiteness property means that each 
bounded subset of $\RR^\ell$ is intersected by only finitely 
many cells. So, the set of cells of a tessellation $T$ 
is necessarily countably infinite. 
We put $C\in T$ for a cell $C$ of the tessellation $T$. 
A tessellation can as well be considered as the closed 
subset $\partial T=\bigcup_{\{C\in T\}} \partial C$ which is
the union of the cell boundaries. There is an obvious one-to-one
relation between both ways of description of a tessellation, and their
measurable structures can be related appropriately, see 
\cite{martnag, schn/weil}. We denote by $\T$ the set of all 
tessellations of $\RR^\ell$.

\medskip

Let ${\cal C}$ be the set of all compact subsets of $\RR^\ell$.
We endow $\T$ with the Borel $\sigma$-algebra ${\cal B}(\T)$ of the
Fell topology, namely
$$
{\cal B}(\T)=\sigma \left( \{ \{ T\in
\T :\, \partial T\cap A=\emptyset \} :\, A\in {\cal C} \} \right)\, .
$$
(As usual, for a class of sets $\I$ we denote by $\sigma({\I})$ the smallest
$\sigma$-algebra containing $\I$.) 

\medskip

Let $\F$ be the family of closed sets of $\RR^\ell$. When $\F$ is endowed
with the Fell topology (for definitions and properties see \cite{schn/weil}),
it is a compact Hausdorff space with a countable base, so it is metrizable.
Also the class of nonempty closed sets
$\F'=\F\setminus \{\emptyset\}$ endowed with the restricted Fell topology
is a complete separable metric space, so for any nonatomic probability measure
$P$ on $(\F',{\cal B}(\F'))$, the completed probability space 
$(\F',{\cal B}(\F'), P)$ is Lebesgue. 
Each tessellation $T\in \T$, as a countable collection of polytopes
is a closed set in $\F'$. Furthermore in
Lemma 10.1.2. in \cite{schn/weil} it was shown that $\T\in {\cal B}(\F')$,
so for any nonatomic probability measure
$P$ on $(\T, {\cal B}(\T))$, the completed probability space
$(\T, {\cal B}(\T), P)$ is Lebesgue. (For more detailed arguments 
see \cite{martnag}, Section 1.3.)

\medskip

Let $W$ be a window in $\RR^\ell$.
The tessellations of $W$ are defined similarly and
the class of them is denoted by $\T\wedge W$. If 
$T\in \T$ we denote by $T\wedge W=\{C\cap W: C\in T\}$ 
the induced tessellation on $W$. 
The tessellation $T\wedge W$ has a finite number of cells 
because $T$ is locally finite and we put
$$
\#(T\wedge W): \hbox{ number of cells of } T\wedge W.
$$
The boundary of $T\wedge W$ is  
$\partial (T\wedge W)= (\partial T \cap W)\cup {\partial W}$.
We introduce the following $\sigma$-algebra,
$$
{\cal B}(\T\wedge W)=\sigma \left( \{ \{ T\in \T\wedge W :\, 
\partial T\cap A=\emptyset \} :\, A\subseteq W,\, A\in {\cal C} \} \right) .
$$
Also, for any nonatomic probability measure
$P$ on $(\T\wedge W, {\cal B}(\T\wedge W))$, the completed probability space
$(\T\wedge W, {\cal B}(\T\wedge W), P)$ is Lebesgue.

\medskip

We note that for another window  $W'\subseteq W$
we have $T\wedge W'=(T\wedge W)\wedge W'$.

\section{Hyperplanes}
\label{section3}

Let $\hH$ be the set of all hyperplanes in $\RR^\ell$, we will
define a parameterization of it.
Let $\|\cdot \|$ be the Euclidean norm, $\langle \cdot , \cdot \rangle$
be the inner product, $\RR_+ =[0,\infty)$ and 
${\SaS}_+^{\ell -1}=\{ x\in \RR^\ell: \| x\| =1\}\cap (\RR^{\ell -1}\times \RR_+)$
be the upper half unit hypersphere in $\RR^\ell$. Define
\begin{equation}
\label{paramx1}
H(\alpha ,u )= \{ x\in \RR^\ell: \langle x,u \rangle =\alpha \},
\quad  \alpha \in \RR , u\in {\SaS}_+^{\ell -1},
\end{equation}
which is the hyperplane with normal direction $u$
and signed distance (in direction $u$) $\alpha$ from the origin.
 Thus we can write
\begin{equation}
\label{paramx2}
\hH = \left\{ H(\alpha ,u ):\,
(\alpha ,u)\in \RR \times {\SaS}_+^{\ell -1} \right\} 
\end{equation}
and on $\hH$ we use the $\sigma$-algebra that is induced
from the Borel $\sigma$-algebra on the parameter space.
Any hyperplane generates two closed half-spaces
$$
H^-(\alpha ,u )= \{ x\in \RR^\ell: \langle x,u \rangle \leq \alpha \}
\; \hbox{ and } \;
H^+(\alpha ,u )= \{ x\in \RR^\ell: \langle x,u \rangle \geq\alpha \} .
$$
For an hyperplane $H$ the above notions are written by $H^-$ and $H^+$ for short.
We define
$$
[B]=\{H\in \hH: \, H\cap  B\neq \emptyset \} \hbox{ for } 
\, B\in {\cal B}(\RR^\ell).
$$

Now, let $\Lambda$ be a (non-zero) measure
on the space of hyperplanes $\hH$. 

\medskip

\subsection{Assumptions on $\Lambda$}
We assume:

\begin{itemize}
\item[(i)] $\Lambda$ is translation invariant;

\item[(ii)] $\Lambda$ possesses the following locally finiteness property:
\begin{equation}
\label{locfin1}
\Lambda([B])\!<\!\infty , \,  \hbox{ for all bounded sets }
B\in {\cal B}(\RR^\ell)\,; 
\end{equation}

\item[(iii)] the support of $\Lambda$ is such that there is
no line in $\RR^\ell$ with the property that all the hyperplanes
of the support are parallel to it.

\end{itemize}

The image of a non-zero, locally finite and
translation invariant measure $\Lambda$ with respect to
the parameterization (\ref{paramx1}), (\ref{paramx2}), 
can be written as the product measure
\begin{equation}
\label{prodmeas}
\gamma \cdot \lambda \otimes \theta,
\end{equation}
where $\gamma >0$ is a constant, $\lambda$ is 
the Lebesgue measure on
$\RR$ and $\theta$ is a  probability measure on
${\SaS}_+^{\ell -1}$ (cf, e.g. \cite{schn/weil}, Theorem 4.4.1 and Theorem 13.2.12).

\medskip

From the properties of $\Lambda$ there is no one-dimensional
subspace $L_1$ of $\RR^\ell$ such that the support of $\theta$ equals
$ L_1^\bot \cap {\SaS}_+^{\ell -1}$ (where $L_1^\bot$ denotes the 
orthogonal complement of $L_1$). This property
allows to obtain a.s. bounded cells in STIT tessellations, cf.
\cite{schn/weil}, Theorem 10.3.2, which can also be applied to STIT.

\medskip

From (\ref{locfin1}) we find that the space $(\hH, {\cal B}(\hH), \Lambda)$ 
is $\sigma-$finite. In fact,
for an increasing sequence of widows $(W_n: n\in \NN=\{0,1,..\})$
covering $\RR^\ell$ (that is $\RR^\ell=\bigcup_{n\in \NN} W_n$)
we have $\hH=\bigcup_{n\in \NN} [W_n]$ and $\Lambda([W_n])<\infty$
for all $n\in \NN$.

\medskip

Let $W$ be a window. Since $\Int W\neq \emptyset$ we get $\Lambda([W])>0$. 
Then, $0<\Lambda([W])<\infty$ and we can define
$$
\hLambda_{[W]}=\Lambda([W])^{-1}\Lambda_{[W]}
$$
the (normalized) probability measure associated to $\Lambda_{[W]}$, 
the restriction of $\Lambda$ to $[W]$.
Since $\Lambda$ is translation invariant we have that $\hLambda_{[W]}$ 
is non-atomic, see \cite{martnag}. Hence, regarding the properties of 
the parameter space, 
which are inherited by the space of hyperplanes we have that 
\begin{equation}
\label{pLebs}
([W],{\cal B},\hLambda_{[W]}) \hbox{ is a Lebesgue probability space }.
\end{equation}

For $T\wedge W\in \T \wedge W$ we define
\begin{equation}
\label{mass1}
\zeta(T\wedge W) = \sum_{C\in T\wedge W} \Lambda ([C]).
\end{equation}

\section{A tree structure}
\label{section4}

Let us introduce the rooted binary trees. First we 
set $\NN^*=\{1,2,...\}$.
Let $\eE=\{-,+\}$ be a two symbol alphabet
and for $k\in \NN^*$ let $\eE^k$ be the set of sequences
(or words) of length $k$. They describe paths in the tree starting from the 
root. We take $\eE^0=\{o\}$ a singleton where $o$
is the empty word. We define $\eE^*=\bigcup_{k\in \NN} \eE^k$.
By $\ve$ we denote an element of $\eE^*$ and we say it has level $k$ if
$\ve\in \eE^k$.

\medskip

Let $\ve=(e_1,...,e_k)\in \eE^k$ for $k\in \NN$ (so $\ve_0=(o)$).
The successors of $\ve$ are the two elements in
$\hbox{Succ}(\ve)=\{(e_1,...,e_k, e_{k+1}): e_{k+1}\in \eE\}$
and $\ve$ is called the predecessor of each of its successors.
For $\ve\in \eE^*\setminus \{o\}$ we denote by Pred($\ve$)
its predecessor. Note that $\hbox{Succ}(o)=\eE$.

\medskip

For simplicity, we will often omit brackets and commas and write 
$\ve=e_1...e_k$ for $\ve=(e_1,...,e_k)$.

\medskip

It is useful to have a total order  $\le$ on $\eE^*$
compatible with the levels.
We fix $(\eE^*,\le)$ as the totally ordered set that satisfies:
\begin{eqnarray}
\label{lexorder}
&{}& \ve \in \eE^j,\ \ve\,' \in \eE^k,\   j<k \, \Rightarrow\,  
\;\, \ve < \ve\,';\\
\nonumber
&{}&  \hbox{ on } \eE : - \, < \, + \,; 
\, \hbox{ inducing } \le \hbox{ on }
 \eE^k, \, k\in \NN^*\,, \hbox{ the lexicographical order}.
\end{eqnarray}
In particular the empty word is the minimal one:
$o<\ve$ for all $\ve\in \eE^*$ with $\ve\neq o$; and
$\ve \,-< \ve \,+$ is the order between the successors of $\ve$.

\medskip

Now we put all paths of a binary tree into ordered tuples $R$. These tuples 
will later be useful for the description of STIT processes.
Note that there can be different such tuples referring to the same tree. 
These tuples can be interpreted as a protocol in which order the  edges of 
the tree 'grow'. It is 
assumed that the pairs of edges to the two successors of a node appear 
simultaneously, but different pairs cannot 'grow' simultaneously.

\medskip

Below for a finite sequence $R=(r_0,...,r_{2k})\in (\eE^*)^{2k+1}$
we denote by $\{R\}=\{r_i: i=0,...,2k\}$ the set of its components and
by $|R|$ the cardinality of $\{R\}$. So, $|R|=2k+1$ means
that the values in $R$ are pairwise different.

\medskip

We shall define the following set $\R$:
\begin{eqnarray}
\nonumber
&{}& \R=\bigcup_{k\in \NN} \R_k \hbox{ with } \R_k\subset (\eE^*)^{2k+1}
 \hbox{ and } R=(r_0,...,r_{2k})\in (\eE^*)^{2k+1}
\hbox{ satisfies } \\
\label{defset}
&{}& R\in \R_k \Leftrightarrow
\, r_0=o\,,\; |R|=2k+1, \hbox{ and } \\
\nonumber
&{}& 
\forall k\ge 1\,   \forall l\in \{1,...,k\} \, \exists j_l\in \{0,...,2l\!-\!2\} :\, 
 \{r_{2l-1},r_{2l}\}=\hbox{Succ}(r_{j_l}), \, r_{2l-1}<r_{2l}.
\end{eqnarray}
Note that
$$
\Big[\{r_{2l-1},r_{2l}\}=\hbox{Succ}(r_{j_l}), \, r_{2l-1}<r_{2l}\Big] 
\, \Leftrightarrow \, \Big[(r_{2l-1},r_{2l})=(r_{j_l}\,-,r_{j_l}\,+)\Big].
$$
For $R\in \R$ and $\ve\in \eE^*$ we have that:
$$
\hbox{Succ}(\ve)\cap \{R\}=\emptyset \, \hbox{ or }
\, \hbox{Succ}(\ve)\cap \{R\}=\hbox{Succ}(\ve).
$$
We define the set of leaves of $R$ by
\begin{equation}
\label{defleav}
\Le(R)=\{r\in \{R\}: \, \hbox{Succ}(r)\cap \{R\}= \emptyset\}.
\end{equation}
This is the set of elements in $\{R\}$ such that both
successors are not in $\{R\}$. 
When $R=(o)$ we have $\Le(R)=\{o\}$. For $k>0$,
and $R=(r_0,...,r_{2k})\in \R_k$ we have in particular
$$
\{r_{2k-1}, r_{2k}\}\subseteq \Le(R).
$$
If $|R|=3$ we have $R=(o,-,+)$ and $\Le(R)=\{-,+\}$.
For $|R|=5$ we have $R=(o,-,+, --, -+)$ or $R=(o,-,+, +-, ++)$,
in the first case $\Le(R)=\{+,--,-+\}$ and in the second one
$\Le(R)=\{-,+-,++\}$.

\medskip

Take $R=(r_0,...,r_{2k})\in \R_k$. For all $s\in \{0,...,k\}$
we define $R^{(s)}=(r_0,...,r_{2s})\in \R_s$. E.g., if 
$R=(o,-, +)$ we have $R^{(0)}=(o)$ and $R^{(1)}=R$.

\medskip

For $k>0$ and $s\in \{0,...,k-1\}$ there exists $j_{s+1}\le 2s$
such that $\hbox{Succ}(r_{j_{s+1}})=\{r_{2s+1},r_{2(s+1)}\}$. 
We denote $r^*_s=r_{j_{s+1}}$ and then we have
\begin{equation}
\label{hoj1}
\forall k>0, \,s\in \{0,...,k-1\}:\;\;
\Le(R^{(s+1)})=(\Le(R^{(s)})\setminus \{r^*_s\})
\cup \{r_{2s+1}, r_{2s+2}\},
\end{equation}
that is $r^*_s\in \Le(R^{(s)})$ is substituted by its successors.
Then, 
\begin{equation}
\label{countxx}
\forall\, R=(r_0,...,r_{2k})\in \R_k: |\Le(R)|=k+1
\hbox{ and } |\Le(R^{(s)})|=s+1 \hbox{ for } s=0,...,k.
\end{equation}
Note that we always have $r^*_0=o$.

\section{Construction of STIT tessellations in a window: main properties}
\label{section5}

A STIT tessellation is defined as a homogeneous (i.e. spatially stationary) 
tessellation with a distribution that is invariant under rescaled iteration 
(or nesting)
of tessellations. A precise definition was given in \cite{nagel/weiss05} where also the
existence was shown (by construction) as well as 
the uniqueness of its law if a hyperplane 
measure $\Lambda$ is given. Meanwhile, several equivalent constructions 
of STIT tessellations in a bounded window
are published. Here we start with one of these constructions.

\medskip

On every window $W$ and for every hyperplane measure $\Lambda$ 
satisfying the assumptions (i), (ii), (iii), formulated 
in Section \ref{section3},  there is a STIT tessellation process
$Y\wedge W =((Y \wedge W)_t: t\ge 0)$ associated to
$\Lambda_{[W]}$, that is now constructed.

\medskip

Let us take two independent families of independent random variables
$(Z(\ve):\ve\in \eE^*)$ and $(G_n: n\in \NN^*)$, 
with $Z(\ve)\sim \,$ Exponential$(1)$
and $G_n\sim$ \,$\hLambda_{[W]}$.
So $\PP(Z(\ve)>t)=e^{-t}$ for all $t\ge 0$. 
We note that $\lambda^{-1}Z(\ve)\sim \,$ Exponential$(\lambda)$ for $\lambda >0$. 

\medskip

Now we define cells $C(\ve)$  which are later used to describe the states 
of the STIT tessellation process.

\begin{enumerate}

\item[Step I:] $C(o)=W$. 

\item[Step II:] Define
\begin{eqnarray}
\label{kappax}
&{}& \forall \ve\in \eE^*:\;\; H(\ve)=G_{\kappa(\ve)} \hbox{ where }
\, \kappa(o)=1 \, \hbox{ and for }  \ve\neq o: \\
\nonumber
&{}& \kappa(\ve)=\inf\left\{{n}: G_n\in [C(\ve)],\,
n>\max\{\kappa({\ve\,'}): {\ve\,'}<\ve\}\right\}.
\end{eqnarray}
So $(H(\ve):\ve\in \eE^*)$ and $(\kappa(\ve): \ve\in \eE^*)$
are well-defined a.s. 

\item[Step III:] For $\ve\in \eE^*$, define
$C(\ve \,-)=C(\ve)\cap H^-(\ve)$, $C(\ve \,+)=C(\ve)\cap H^+(\ve)$.

\item[Step IV:] $C(o)$ is born at time $t_b(o)=0$, its
lifetime is $t_l(o)\sim \Lambda ([C(o)])^{-1}\, Z(o)$
and so dies at $t_l(o)$, i.e. at that time it is divided by $H(o)=G_1$. 
A cell $C(\ve)$, with $\ve\in \eE^*$, $\ve\neq o$, is born 
at time $t_b(\ve)=t_b(\hbox{Pred}(\ve))+t_l(\hbox{Pred}(\ve))$,
has lifetime $t_l(\ve)\sim \Lambda ([C(\ve)])^{-1}\,
Z(\ve)$ and dies at time $t_b(\ve)+t_l(\ve)$. At that time it is divided 
by $H(\ve)=G_{\kappa(\ve)}$ into $C(\ve \,-)$ and $C(\ve \,+)$.

\end{enumerate}

In step II a rejection method is applied where random hyperplanes
are thrown onto the window until the first time a hyperplane hits $C(\ve)$.
Note that in (\ref{kappax}), the sequence $(\kappa(\ve): \ve\in \eE^*)$ 
is increasing and $(H(\ve):\ve\in \eE^*)$  is an independent family 
conditioned to $H(\ve)\cap [C(\ve)]\neq \emptyset$. Moreover,
\begin{equation}
\label{disHx}
\forall \, \ve\in \eE^*:\;\; H(\ve)\sim \hLambda_{[C(\ve)]}.
\end{equation}

\begin{remark}
\label{rem1}
Let $\ve\neq o$. In the sequence 
$\left(G_n: n> \max \{\kappa(\ve \,'): {\ve\,'}<\ve\}\right)$ 
of independent identically distributed
random hyperplanes with common law $\hLambda_{[W]}$, the first
of these hyperplanes which intersects $[C(\ve)]$ is distributed as
$\hLambda_{[C(\ve)]}$. The random time of attending such an
hyperplane depends on the inverse of the $\Lambda$-measure of $[C(\ve)]$, 
which depends on the size but also on the shape of the cell $C(\ve)$.
\end{remark}

It is easy to see that at any time a.s. at most one cell
dies and so a.s. at most only two cells are born.

\medskip

At each time $t\ge 0$ we define $(Y\wedge W)_t$ as the class
of cells $C(\ve)$ which are alive at time $t$, that is 
$$
(Y\wedge W)_t=\{C(\ve): \ve\in E_t\} \hbox{ where }
E_t=\{\ve\in \eE^*: t_b(\ve)\le t < t_b(\ve)+t_l(\ve) \}.
$$
Since $C(\ve)=C(\ve \,-)\cup C(\ve \,+)$
and $\emptyset=\Int C(\ve \,-)\cap \Int C(\ve \,+)$, it
is easy to see that $(Y\wedge W)_t$ is a tessellation of $W$.
On the other hand, it can be checked that $E_t=\Le(R)$, 
the set of leaves of
some $R\in \R$ defined in Section \ref{section4}.
Note that such an $R$ is not necessarily unique.

\medskip

Let $t\ge 0$, $s>0$. Let us show that $(Y\wedge W)_{t+s}$ is conditionally 
independent from $((Y\wedge W)_{v}, v<t)$, conditioned on $(Y\wedge W)_{t}$.
From definition 
$$
(Y\wedge W)_{t+s}\cap \left\{ C(\ve) :\, 
\ve\in \left( \bigcup_{v<t} E_v \right) \setminus E_t \right\}  =\emptyset .
$$
On the other hand 
\begin{eqnarray*}
&{}& (Y\wedge W)_{t+s}\subseteq \bigcup_{\ve\in E_t}
\{C({\ve\, '}): {\ve\, '}\in \hbox{Succ}^*(\ve)\} \, \hbox{ where }\\
&{}&
\hbox{Succ}^*(\ve)=\{{\ve\, '}: \exists k\ge 1, \exists 
{\ve\, ^1},...,{\ve\, ^{k-1}}, \forall \, j\!=\!1,...,k:
{\ve\, ^j}\in \hbox{Succ}(\ve\, ^{j-1}),
\, {\ve\, ^0}=\ve, {\ve\, ^k}={\ve\, '} \}.
\end{eqnarray*}
Then, the memoryless property of the exponential distribution,
implies that $Y\wedge W$ is a Markov processes,
$$
\PP\left((Y\wedge W)_{t+s}\in \cdot \, | \, (Y\wedge W)_v, v\in [0,t]\right) 
=\PP\left((Y\wedge W)_{t+s} \in \cdot \, | \, (Y\wedge W)_t \right).
$$
From (\ref{mass1}),
\begin{equation}
\label{defzeta}
\zeta((Y\wedge W)_t) = \sum_{C\in (Y\wedge W)_t} \Lambda ([C]).
\end{equation}
Let 
$$
\tau_t=\inf\{s>0: (Y\wedge W)_{t+s}\neq (Y\wedge W)_t\}.
$$
be the holding time at $t$.
The memoryless property of the exponential distribution 
(and the property that the minimum of finitely many independent exponentially 
distributed random variables is again exponentially distributed 
where its parameter is the sum 
of the parameters of the variables) implies,
\begin{equation}
\label{laexp1}
\tau_t \, | \, \sigma((Y\wedge W)_t) \sim \,
\hbox{ Exponential}(\zeta((Y\wedge W)_t)). 
\end{equation}
At time $t+\tau_t$ a (a.s.) unique cell $C^*_t\in (Y\wedge W)_t$ dies,
we put $C^*_t=C(\ve\,^*)$ with ${\ve\,^*}\in E_t$ a random index. 
Then, two new cells 
$\{C^*_t\cap H^-({\ve\,^*}), C^*_t\cap H^+({\ve\,^*})\}$ 
are born at this time, so
\begin{equation}
\label{uni0}
(Y\wedge W)_{t+\tau_t}=\{C: C\in (Y\wedge W)_t, C\neq C^*_t\}\cup 
\{C^*_t\cap H^-({\ve\,^*}), C^*_t\cap H^+({\ve\,^*})\}.
\end{equation}
Hence,
\begin{equation}
\label{uni1}
Y_{t+\tau_t}  \hbox{ is uniquely defined from } \big[(Y\wedge W)_t,\,
C^*_t\in (Y\wedge W)_t,\, H({\ve\,^*}) \big].
\end{equation}

To any realization of the process $((Y\wedge W)_t: 0\le t\le t_0)$ 
we associate a sequence $(R^{s}\in \R: s=0,..,s_0)$ as follows. 
Define $R^{(0)}=(o)$. Now let
$\tau^0=0$ and for an integer $s\ge 0$ define
$$
\tau^{s+1}=\inf \{t> 0: \tau^s+t\le t_0,
(Y\wedge W)_{\tau^s+t}\neq (Y\wedge W)_{\tau^s}\},
$$
where as usual $\infty=\inf \emptyset$. Then a.s. there exists  
$s_0=\sup\{s: \tau^s\le t_0\}$. 
Then, we define $R=R^{(s_0)}$ by induction as follows: 
for all $0\le s <s_0$ define, 
\begin{equation}
\label{uniR}
R^{(s+1)}=(r_0,...,r_{2s},r^* -   , r^* +) ,
\end{equation}
where $r^* \in \{ r_0,...,r_{2s} \}$ is the unique element in $R^{(s)}$,
such that  $C^*_{\tau^s}=C(r^*)$ according to (\ref{uni0}).

\medskip

In next result we characterize the Markov process $Y\wedge W$
by supplying its holding times, the jump rates and their conditional 
independence. We follow Section 1.1.1. in \cite{ber}. 

\medskip

\begin{proposition}
\label{lawstit}
(I) $Y\wedge W$ is a pure jump Markov process and satisfies
\begin{eqnarray}
\label{lawX0}
&{}& \forall C(\ve)\in (Y\wedge W)_t,\forall H\in [C(\ve)], s>0: \\
\nonumber
&{}& \PP(H(\ve)\in dH, C^*_{t}=C(\ve), \tau_t\in ds \, | \,
(Y\wedge W)_t)=\Lambda_{[C(\ve)]}(dH)
e^{-\zeta((Y\wedge W)_t)\, s}ds.
\end{eqnarray}
(II) We have the consistency property, namely
\begin{equation}
\label{consx}
\forall \hbox{ windows } V\subseteq W:\;\;
(Y\wedge W)\wedge V\sim Y\wedge V \hbox{ where }
(Y\wedge W)\wedge V=((Y\wedge W)_t\wedge V: t\ge 0).
\end{equation}
(III) Moreover, $Y\wedge W$ satisfies the following regeneration property:
For all fixed $t_0\ge 0$ and $C\in Y_{t_0}$
the processes $(Y\wedge W)\wedge C=((Y\wedge W)_t\wedge C: t\ge t_0)$, satisfy
\begin{eqnarray}
\nonumber
&{}&\left((Y\wedge W)\wedge C: C\in Y_{t_0}\right)
\hbox{ are conditionally independent, given } Y_{t_0} \hbox{ and}\\
\label{fragnewx}
&{}&  
(Y\wedge W)\wedge C \, \hbox { is a STIT process on } C, 
\hbox{ associated to } \Lambda .
\end{eqnarray}
\end{proposition}

{\em Proof:}
From (\ref{disHx}) and Step III, we have that
$$
\forall C(\ve)\in (Y\wedge W)_t,\, \forall H\in [C(\ve)]:\;\,
\PP(H(\ve)\in dH  \, | \, C^*_{t}=C(\ve), (Y\wedge W)_t)=\hLambda_{[C(\ve)]}(dH).
$$
Hence, (\ref{lawX0}) will follow once we prove 
$$
\forall C\in (Y\wedge W)_t:
\;\; \PP(C^*_{t}=C, \tau_t\in ds \, | \,
(Y\wedge W)_t)=\Lambda([C])\,
e^{-\zeta((Y\wedge W)_t)\, s}ds\,.
$$
The proof of this last relation is based upon the following fact applied to 
the lifetime variables $Z(\ve)$ of $C(\ve)\in (Y\wedge W)_t$.   
Let $(Z_i: j=1,...,k)$ be independent random variables with 
$Z_i\sim \,$ Exponential$(q_i)$ and $Z=\min\{Z_i: i=1,...,k\}$. 
Then for all $z>0$,
$$
\PP(Z_j\!=\!Z, \, Z\!\in \!dz)=\PP(Z_i\!>\!z, i\!\neq \!j, Z_j\!\in \!dz)
=\frac{q_j}{\sum_{i=1}^k q_i}\PP(Z\in dz)=
q_j \, e^{-(\sum_{i=1}^k q_i)z} dz. 
$$

The consistency property (\ref{consx}) was already shown in \cite{nagel/weiss05}.    

\medskip

The proof of the regeneration property (\ref{fragnewx}) 
follows from the consistency and the memoryless property of the 
exponential distribution
applied to $C(\ve\, ')\in Y_{t_0}\wedge W$, which is
$$
\forall t\ge t_0:\; \;
\PP(t_b(\ve)+t_l(\ve)> t \, | \, t_b(\ve) \le t_0 < t_b(\ve)+t_l(\ve) \}=
e^{-\Lambda([C(\ve)])(t-t_0)}.
$$
For $t\ge t_0$ and $C(\ve\, ')\in {Y_{t_0}\wedge W}$, 
$(Y\wedge W)\wedge C(\ve\, ')=((Y\wedge W)_t\wedge C(\ve\, '): t\ge t_0)$ 
satisfies
$$
(Y\wedge W)_t\wedge C(\ve\, ')=\{C(\ve): \ve\in {Succ}^*(\ve\, '),
t_b(\ve)\le t_0 \le t < t_b(\ve)+t_l(\ve) \}.
$$

From the memoryless property we find that 
$(Y\wedge W)\wedge C(\ve\, ')$ is a STIT process. 
Finally, the interiors of the cells 
$C\in \bigcup_{t\ge 0}(Y\wedge W)_t\wedge C(\ve\, ')$ are contained in 
$\Int C(\ve\, ')$, and so they are pairwise disjoint  
as $C(\ve\, ')$ varies in $Y_{t_0}\wedge W$. We deduce that the processes 
$(Y\wedge W)\wedge C(\ve\, ')$, $C(\ve\, ')\in Y_{t_0}\wedge W$,
are conditionally independent. \hfill $\Box$

\medskip

Note that (\ref{lawX0}) implies,
$\tau_t \, | \, (Y\wedge W)_t  \, \sim \, \hbox{Exponential}(\zeta((Y\wedge W)_t))$,
and that for all $C(\ve)\in (Y\wedge W)_t$, $H\in [C(\ve)]$  holds
\begin{equation}
\label{fragnewselectdivi}
\PP(H(\ve)\in dH, C^*_{t}=C(\ve)  \, | \, (Y\wedge W)_t)=
\hLambda_{[C(\ve)]}(dH)\, \frac{\Lambda([C(\ve)])}{\zeta((Y\wedge W)_t)},
\end{equation}
in particular
$\PP(C^*_{t}=C  \, | \, (Y\wedge W)_t)={\Lambda([C])}/{\zeta((Y\wedge W)_t)}$.

\medskip

In \cite{nagel/weiss05} it was shown that the STIT process $Y\wedge W$ has
no explosion. In fact the process of number of cells
$\#\{C: C\in (Y\wedge W)_t \}$ is
stochastically dominated by a birth chain
$(M(t): t\ge 0)$ starting from $M(0)=1$
with linear birth rates $b_n=n\, \Lambda([W])$. Since $(M(t): t\ge 0)$
does not explode we deduce that the process of the number of cells 
does not explode too.
This also follows straightforwardly from Lemma 1.1 in \cite{ber}.

\medskip

The consistency property (\ref{consx}) implies the existence of
a probability measure on $\T^{\RR_+}$ endowed with the product $\sigma-$field,
and it defines the distribution of a process $Y$, which is Markov 
and it satisfies $Y_t\wedge W\sim (Y\wedge W)_t$ for all 
$t\in \RR_+$ and all windows $W$.
(See  \cite{nagel/weiss05}).

\medskip

A global construction for a STIT process $Y$
was provided in \cite{m/n/w08a}.

\begin{proposition}
\label{propglobal}
The process $Y=(Y_t: t>0)$ is a Markov process which is the STIT tessellation 
process on $\RR^\ell$ associated to $\Lambda$. Its marginals $Y_t$ 
take values in $\T$ and $(Y_t\wedge W: t\ge 0)\sim ((Y\wedge W)_t: t\ge 0)$.

\medskip

Moreover $Y$ satisfies the regeneration property: For all $t_0> 0$ and 
$C\in Y_{t_0}$ the processes $Y\wedge C=(Y_t\wedge C: t\ge t_0)$,
 satisfy
\begin{eqnarray}
\label{fragnewxgen}
\left(Y\wedge C: C\!\in \!Y_{t_0}\right)
\hbox{ are conditionally independent given } Y_{t_0}, \hbox{ and } \nonumber\\
Y\wedge C \, \hbox { is a STIT process on } C \hbox{ associated to } \Lambda .
\end{eqnarray}
\end{proposition}

The proof of the regeneration property (\ref{fragnewxgen}) is 
straightforward from (\ref{fragnewx}). 
We notice that this regeneration property (\ref{fragnewxgen}) 
is equivalent to the stable-under-iteration property. For this 
last property see 
(\cite{martnag, st}). On the other hand (\ref{fragnewxgen}) is, 
once written appropriately, the branching property of a 
fragmentation chain described in Proposition 1.2 $(i)$ \cite{ber}.

\section{The marginal distribution of STIT tessellations in a window}
\label{section6}

In (\ref{uniR}) we  used the tuples $R$ defined for a rooted binary tree
in Section \ref{section4} to index the sequences of tessellations associated
constructed in Steps I and III in Section \ref{section5}.

\medskip

For each $R=(r_0,...,r_{2k})\in \R_k$ 
we have defined a sequence $R^{(s)}=(r_0,...,r_{2s})\in \R_s$ for $s=0,...,k$.
We will associate to each $R=(r_0,...,r_{2k})\in \R_k$ a 
sequence of tessellations of $W$, 
denoted by $T(R)=(T(R^{(s)}): s=0,...,k)$.
We will do it by describing the family of cells
of each $T^{(s)}$, by using induction on $s=0,...,k$. 

\medskip

We define
$$
T^{(s)}(R)=\{C(r): r\in \Le(R^{(s)})\}.
$$
Note that $T^{(0)}(R)=\{C(o)\}$ because $R^{(0)}=(o)$.
Let $k\ge 1$.
From (\ref{hoj1}) we get for $s=0,...,k-1$:
$$
T^{(s+1)}(R)=\{C\in T^{(s)}(R)\setminus C(r^*_s)\}
\cup \{C(r_{2s+1}), C(r_{2s+2})\}.
$$
That is the tessellation $T^{(s+1)}(R)$ results from dividing the
cell $C(r^*_s)\in T^{(s)}(R)$ into the two cells 
associated to its successors (here, the dividing hyperplane $H(r^*_s)$, 
see (\ref{kappax}), is not indicated in the notion).

\medskip

We have $\#(T^{(s)}(R))=s+1$, see (\ref{countxx}), and
$\zeta(T^{(s)}(R))=\sum_{C\in T^{(s)}(R)} \Lambda([C])$, see (\ref{defzeta}).

\medskip

Let $A\in {\cal B}(\T)$. The STIT process
$Y\wedge W=((Y\wedge W)_t: t\ge 0)$ satisfies
$$
\PP((Y\wedge W)_t\in A)= \PP((Y\wedge W)_t=\{W\}, \{W\}\in A)+
\sum_{k\in \NN^*}\PP(\#(Y\wedge W)_t=k+1, (Y\wedge W)_t\in A).
$$
To describe the summands for $k\in \NN^*$, 
we must take into account that at any time $s$, the tessellation 
$(Y\wedge W)_s$ attends a random time $\tau_s$ for the
division of one of its cells, and this time satisfies
$\tau_s \, | \, (Y\wedge W)_s  \, \sim \, \hbox{Exponential}(\zeta((Y\wedge W)_s))$
 and the cell $C^*_s$ of $(Y\wedge W)_s$ divided at 
time $s+\tau_s$ is chosen by 
$\PP(C^*_{t}=C  \, | \, (Y\wedge W)_s)={\Lambda([C])}/{\zeta((Y\wedge W)_s)}$. 
Therefore, using (\ref{lawX0}), we obtain
\begin{proposition} For $k\in \NN^*$
\begin{eqnarray}
\label{margfor1}
&{}&\PP(\#(Y\!\wedge \!W)_t\!=\!k\!+\!1, (Y\!\wedge \!W)_t \!\in \!A)\\
\nonumber
&{}&= \sum_{R\in \R_k}\! \int\!\!\! d\hLambda_{[C(r^*_0)]}(H_1) 
\int\!\!\! d\hLambda_{[C(r^*_1)]}(H_2)...\int\!\!\!
d\hLambda_{[C(r^*_{k-1})]}(H_k)
\int_0^t \!\!\!dw_1 \!\!\int_0^{t-w_1}\!\!\!\!\!\!\!\!dw_2...
\int_0^{t-\sum_{s=1}^{k-1}w_j} \!\!\!\!dw_{k} \\
\nonumber
&{}&
\;\;\times \left(\prod_{s=0}^{k-1} \Lambda([C(r^*_s)])\right)
\, exp \left[ -\sum_{s=0}^{k-1}\zeta(T^{(s)}(R)) w_{s+1}\right]
exp \left[ {-\zeta(T^{(k)}(R))(t-\sum_{s=1}^k w_j)} \right]
{\bf 1}_{A}(T^{(k)}(R)).
\end{eqnarray}
\end{proposition}

\section{Revisiting the STIT tessellation process in a window}
\label{section7}

Let us give another equivalent construction of the STIT process $Y \wedge W$
that will be useful  to understand
the construction done in the next Section.

\medskip

Let us consider three independent sequences
$(U_n: n\in \NN^*)$, $(V_n: n\in \NN^*)$, $(G_n: n\in \NN^*)$,
of independent identically  distributed random variables, such that
$U_n\sim \,$Uniform$[0,1)$, $V_n\sim \,$Exponential$(1)$,
$G_n\sim \hLambda_{[W]}$. We start with a construction of a 
sequence $(\Y_n,\ n\in \NN)$ of tessellations in $W$.

\medskip

The algorithm of the construction is:

\medskip

\noindent Step $n=0$: $\Y_0=\{C(o)\}$ with $C(o)=W$. Let 
$\zeta(\Y_0)=\Lambda([C(o)])$ and $\kappa_0=0$.

\medskip

\noindent  Step $n+1$:

Assume $\Y_n=\{C(\ve): \ve\in E_n\}$, $n\ge 0$, has been defined
with $E_n$ a set of the form $E_n=\Le(R)$ 
with some $R\in \R_{n}$, and so $|E_n|=n+1$. 
We also assume $\kappa_n\in \NN^*$ has been defined.
Let
$$
\zeta(\Y_n)=\sum_{C \in \Y_n} \Lambda([C])=
\sum_{\ve \in E_n} \Lambda([C(\ve)]).
$$
Since $E_n$ is totally ordered by (\ref{lexorder}), also the class of
cells $\{C: C\in \Y_{n}\}$ is totally ordered.
We define a partition of $[0,1)$ by
\begin{equation}
\label{new0}
[0,1)=\bigcup_{\ve\in E_{n}}[a^n_\ve, b^n_\ve) \,
\hbox{ with } \, b^n_\ve - a^n_\ve=\zeta(\Y_n)^{-1}\Lambda([C(\ve)]),
\end{equation}
where the intervals $[a^n_\ve, b^n_\ve)$ and $[a^n_{\ve\,'},b^n_{\ve\,'})$
are consecutive when ${\ve\,'}$ is the element following $\ve$ in $E_{n}$
with respect to the total order $\le$. 
We define a random cell $C^*_n\in \Y_n$ by
\begin{equation}
\label{newx2}
\forall {\ve}\in E_n:\;\;
C^*_n=C(\ve) \Leftrightarrow U_{n+1}\in [a^n_{\ve},b^n_{\ve}).
\end{equation}
Hence
\begin{equation}
\label{newx2'}
\forall \ve \in E_n:\;\;
\PP(C^*_n=C(\ve) \, | \, \Y_n)=\zeta(\Y_n)^{-1}\Lambda([C(\ve)]).
\end{equation}
We denote by ${\ve\, ^*}\in E_n$ the random index such that $C({\ve\, ^*})=C^*_n$.
Note that
\begin{equation}
\label{newx1}
C^*_n\in \sigma(\Y_n,U_{n+1})
\end{equation}

We define the random hyperplane $H_{n+1}$ in a similar way as in (\ref{kappax}),
so
\begin{equation}
\label{newx1''}
H_{n+1}=G_{\kappa_{n+1}} \hbox{ where }
{\kappa_{n+1}}=\min\{j>\kappa_n: G_j\in [C^*_n]\}.
\end{equation}
By definition $H_{n+1}\sim \hLambda_{[C^*_n]}$. Note that $H_1=G_1$. 
Obviously $(\kappa_{n}: n\in \NN^*)$ 
is an increasing sequence of random times.
We note that $\kappa_{n+1}$ is a  stopping time with respect to the filtration 
$(\sigma(G_j,\kappa_{n},C^*_n): j\in \NN^*)$.

\medskip

The tessellation $\Y_{n+1}$ is formed from $\Y_n$ by the 
division of the random cell $C^*_n$ of $\Y_n$ by $H_{n+1}$,
giving
\begin{equation}
\label{newx5}
\Y_{n+1}=(\{C\in \Y_{n}\}
\setminus \{C^*_n\})\cup \{C^*_n\cap H_{n+1}^-, C^*_n \cap H_{n+1}^+\}.
\end{equation}
So $\Y_{n+1}$ is indexed by
$E_{n+1}=(E_n\setminus \{\ve\, ^*\})\cup
\hbox{Succ}({\ve\, ^*})\subset \eE^*$, and so
$E_{n+1}=\Le(R_{n+1})$ for some (uniquely determined) $R_{n+1}\in \R_{n+1}$.
This shows $(\Y_n: n\in \NN)$ is well-defined.
Notice that $\Y_{n+1}\in \sigma(\Y_n,U_{n+1},H_{n+1})$
and then by recursion we get
\begin{equation}
\label{newx7}
\Y_{n+1}\in \sigma(U_k,H_k: k\le n+1) .
\end{equation}

Now, use (\ref{newx2'}) and (\ref{newx1''}) 
to get that for all $\K\in {\cal B}([W])$ we have,
\begin{eqnarray}
\nonumber
&{}& \PP(H_{n+1}\in \K, C^*_n=C(\ve) \, | \, \Y_n)
= \PP(H_{n+1}\in \K \, | \, C^*_n=C(\ve), \Y_n)\PP(C^*_n=C(\ve) \, | \, \Y_n)\\
\nonumber
&{}& =\PP(H_{n+1}\in \K\cap [C(\ve)])\, \zeta(\Y_n)^{-1}\Lambda([C(\ve)])
=\frac{\Lambda(\K\cap [C(\ve)])}{\Lambda([C(\ve)])} \zeta(\Y_n)^{-1}
\Lambda([C(\ve)])\\
\label{newx6'}
&{}&=\zeta(\Y_n)^{-1}\Lambda([C(\ve)]\cap \K).
\end{eqnarray}

Since $\sigma(U_n, G_n:n\in \NN) \, \qperp \,  \sigma(V_n:n\in \NN)$,
from (\ref{newx1}), (\ref{newx1''}) and  (\ref{newx7})
the above random objects satisfy the relation
\begin{equation}
\label{newx8}
\sigma(\Y_n, U_n, G_n, H_n: n\in \NN) \, \qperp \,  (V_n:n\in \NN).
\end{equation}

Define the sequence of jump times by
$$
\forall n\in \NN^*:\quad \pi_n=(\zeta(\Y_{n-1}))^{-1}V_n,
$$
which are conditionally distributed as,
\begin{equation}
\label{newx9}
\pi_n \, | \, \Y_{n-1} \sim \hbox{Exponential}(\zeta(\Y_{n-1})).
\end{equation}
Since $(V_n: n\in \NN^*)$ is a sequence of independent random variables,
$(\pi_n: n\in \NN^*)$ is conditionally independent given 
$\sigma(\zeta(\Y_n): n\in \NN)$.
From (\ref{newx8}) and (\ref{newx9}) we get
\begin{equation}
\label{indht1}
(\pi_{n+1} \, \qperp \, \sigma(C^*_{n}, H_{n+1})) \; | \, \Y_n.
\end{equation}

Define the sequence of times
$$
\aS_0=0 \, \hbox{ and } \; \aS_n=\aS_{n-1}+\pi_n=\sum_{j=1}^{n} \pi_j \,
\hbox{ for } n\in \NN^*.
$$
The proof ensuring that
$\lim\limits_{n\to \infty}\aS_n=\infty \,$ $\PP$-a.s.,
is the same as the one where we proved that $Y\wedge W$ has no explosion. 
In fact, $\#(\Y_n)=n+1$ implies $\zeta(\Y_n)\le (n+1)\Lambda([W])$.
Hence, the process $(N(t): t\ge 0)$ given by 
\begin{equation}
\label{defNx}
N(t)=\sup\{n\in \NN: \aS_n\le t\}, 
\end{equation}
is stochastically dominated by a birth chain $(M(t): t\ge 0)$ starting 
from $M(0)=1$ with linear birth rates $b_n=n\, \Lambda([W])$, so 
$(N(t): t\ge 0)$ does not explode.

\medskip

By using the sequences $(\Y_n:n\in \NN)$ and $(\aS_n:n\in \NN)$ 
we define the tessellation process $\tY\wedge W=((\tY\wedge W)_t: t\ge 0)$ by
\begin{equation}
\label{newx11}
({\tY}\wedge W)_t=\Y_{n} \, \hbox{ when } t\in
\left[\aS_{n}, \aS_{n+1}\right), n\in \NN.
\end{equation}
From $\lim\limits_{n\to \infty}\aS_n=\infty\,$ $\PP-$a.s. 
we get that ${\tY}\wedge W$ is well-defined
for all times $t\ge 0\;$ $\,\PP-$a.s. We also have
$({\tY}\wedge W)_{\aS_n}=\Y_n$
for all $n\in \NN$ and so $(\aS_n: n\in \NN^*)$ is the sequence
of times of jumps of $\tY\wedge W$. 

\begin{proposition}
\label{primeas}
The process $\tY\wedge W$ is a STIT process
associated to $\Lambda_{[W]}$.
\end{proposition}

{\em Proof:}
Since (I) in Proposition \ref{lawstit} completely characterizes the law 
of a STIT process associated to $\Lambda$, 
it is sufficient to show those properties.

\medskip

Let us prove $\tY\wedge W$ satisfies the Markov property. 
For all $t\ge 0$ and $s>0$ we have
$$
\PP(({\tY}\wedge W)_{t+s} \, | \, ({\tY}\wedge W)_u, u\le t)=
\PP(({\tY}\wedge W)_{t+s} \, | \, (\tY\wedge W)_t=\Y_{N(t)}, \aS_{N(t)}).
$$
The memoryless property of the exponential distribution implies
$(\aS_n-t: n>N(t))  \,  
\qperp \, \aS_{N(t)}  \, | \, \Y_{N(t)}$, and so 
$({\tY}\wedge W)_{t+s}  \, \qperp \, \aS_{N(t)} \, | \, \Y_{N(t)}$.
Then
$$
\PP(({\tY}\wedge W)_{t+s} \, | \, ({\tY}\wedge W)_u, u\le t)=
\PP(({\tY}\wedge W)_{t+s} \, | \, (\tY\wedge W)_t=\Y_{N(t)}),
$$ 
so the Markov property is satisfied.

\medskip

The process $\tY\wedge W$ is a jump process. Let us compute
the distribution of the holding time 
$\htau_t=\inf\{s>0: ({\tY}\wedge W)_{t+s}\neq
({\tY}\wedge W)_t\}$, Again by  the memoryless property 
of the exponential distribution we get
$$
\PP(\htau_t>s \, | \, \Y_{N(t)})=\PP(\pi_{N(t)+1}>s \, | \, \Y_{N(t)})=
e^{-\zeta(\Y_{N(t)})\, s},
$$
and so 
$\htau_t \, | \, \Y_{N(t)} \sim \hbox{Exponential}(\zeta(\Y_{N(t)}))$.
Now, from (\ref{indht1}) we deduce the conditional independence relation, 
\begin{equation}
\label{indht2}
\htau_t \, \qperp \, \sigma(C^*_{N(t)}, H_{N(t)+1}) \; | \, \Y_{N(t)},
\end{equation}
and so, from (\ref{newx2'}) and (\ref{newx1''}), we get for all 
$C(\ve)\in \Y_{N(t)}$, $H\in [C(\ve)]$ and $s>0\,$: 
$$
\PP(H(\ve)\in dH, C^*_t=C(\ve),\htau_t\in ds \, | \, \Y_{N(t)})=
\hLambda_{[C(\ve)]}(dH)\, \Lambda([C(\ve)]) \, e^{-\zeta(\Y_{N(t)})\, s} ds.
$$
We have proven relation (\ref{lawX0}), so the result follows. \hfill $\Box$

\section{A new construction of STIT tessellations in a window, point processes}
\label{section8}

Fix a window $W$. 

\medskip

Let us consider three independent sequences
$(U_n: n\in \NN^*)$, $(V_n: n\in \NN^*)$, $(G_n: n\in \NN^*)$, 
of independent identically  distributed random variables, such that
$U_n\sim \,$Uniform$[0,1)$, $V_n\sim \,$Exponential$(1)$,
$G_n\sim \hLambda_{[W]}$.

\medskip

We will construct the STIT 
tessellation process $Y\wedge W$ by using these three 
independent sequences. By using $(G_n: n\in \NN)$ and 
regarding the current state of the tessellation process,
 we will construct a sequence 
of random hyperplanes $(H_n: n\in \NN)$ on 
$[W]$ and all the hyperplanes $H_n$ will be 
effectively used in constructing the STIT, contrary 
to the rejection procedure of previous sections where we must 
wait until a random hyperplane cuts a prescribed cell.
The $(U_n: n\in \NN)$ are used to choose the cell to be divided, 
out of the set of cells which are intersected by the $H_n$.

\medskip

The construction will be done in an iterative way.
For $n\in \NN$, $\Y_n$ is a random tessellation 
of $[W]$ and $\Gamma_n$ is a random measure on $([W],{\cal B}([W])$ defined by
$$
\Gamma_n=\sum_{C\in \Y_n}\Lambda_{[C]}.
$$
Note that $\Gamma_n$ is absolutely continuous
with respect to $\Lambda_{[W]}$ and its Radon-Nikodym derivate
\begin{equation}
\label{relcr0}
\xi_n(H)=\frac{d\, \Gamma_n}{d\, \Lambda_{[W]}}(H) \hbox{ satisfies }
\xi_n(H)=\#\{C\in \Y_n: H\in [C]\}.
\end{equation}
This follows from the partition:
\begin{equation}
\label{relcr1}
[W]=\bigcup_{j=1}^{\#(\Y_n)} {\widehat{\K}}_j \,
\hbox{ with } {\widehat{\K}}_j=\{H\in [W]: \#\{C\in \Y_n: H\in [C]\}=j\}.
\end{equation}

So, the random measures $\Gamma_n$ can be described by the functions $\xi_n$,
which belong to $L^1(\Lambda)$. (For an explicit form of the corresponding 
density on the parameter space see further in Section \ref{remsim}.)  
Since this is a metric 
complete separable space we are in the framework of measurability
described in Section \ref{section1}. By definition,
\begin{equation}
\label{reldep1}
\sigma(\Gamma_n)\subseteq\sigma(\Y_n).
\end{equation}
Consider the probability distribution 
$\hGamma_n=(\Gamma_n([W]))^{-1}\Gamma_n = (\zeta (\Y_n))^{-1}\Gamma_n$
on $([W],{\cal B}([W])$. From (\ref{pLebs}) and  
from Section \ref{section1} we know that 
there is a bimeasurable function $\Xi_n:[W]\to [W]$ such that
$\hLambda_{[W]}\circ \Xi_n^{-1}=\hGamma_n$.

\medskip

The algorithm of the construction is:

\medskip

\noindent Step $n=0$: $\Y_0=\{W\}$ and $\Gamma_0=\Lambda_{[W]}$.
So $\xi_0\equiv 1$ and $\Xi_0=\hbox{Identity}_{[W]}$.

\medskip

\noindent  Step $n+1$: Assume $\Y_n$, $n\ge 0$, has been defined.
We take 
\begin{equation}
\label{reldep2}
H_{n+1}=\Xi_{n}\circ G_{n+1}.
\end{equation}
Its conditional distribution satisfies
\begin{equation}
\label{forvol}
H_{n+1} \, | \, \Y_{n} \sim \hGamma_{n}
\, \hbox{ i.e. } \,
\forall \K\in {\cal B}[W]:
\;\; \PP(H_{n+1}\in \K \, | \, \Y_{n})=\hGamma_{n}(\K)=
\frac{\Gamma_{n}(\K)}{\Gamma_{n}([W])}.
\end{equation}

\medskip

The tessellation $\Y_{n+1}$ is formed from
$\Y_n$ by the division of the random cell $C^*_n$ of $\Y_n$, 
chosen with the help of $H_{n+1}$ and $U_{n+1}$. All the cells of $\Y_n$
which are hit by $H_{n+1}$ have the same probability $\xi_n(H_{n+1})^{-1}$ 
to be chosen for division. Formally, divide the unit interval $[0,1)$ 
into  $\xi_n(H_{n+1})$ intervals of equal
length, namely

\begin{equation}
\label{relcr2}
[0,1)=\bigcup_{\stackrel {\ve\in E_{n}}{H_{n+1}\in [C(\ve )]}}[a^n_\ve, b^n_\ve) \,
\hbox{ with } \, b^n_\ve-a^n_\ve=
\xi_n(H_{n+1})^{-1}
\end{equation}
if $\Y_n = \{ C(\ve ):\, \ve \in E_n \}$ where $E_n=\Le(R_n)$ is 
a set of leaves, for some $R_n\in \R_{n+1}$.
As in (\ref{new0}),
the intervals $[a^n_\ve, b^n_\ve)$ and $[a^n_{\ve\,'},b^n_{\ve\,'})$
are consecutive when ${\ve\,'}$ is the element following $\ve$ in $E_{n}$.
Now, when
\begin{equation}
\label{decx}
U_{n+1}\in [a^n_{\ve}, b^n_{\ve}) 
\end{equation}
we take $C^*_n=C(\ve)$ and divide it by $H_{n+1}$.
Thus $\Y_{n+1}$ is defined as
\begin{equation}
\label{stnplus}
\Y_{n+1}=(\{C\in \Y_{n}\}
\setminus \{C^*_n\})\cup \{C^*_n\cap H_{n+1}^+, C^*_n \cap H_{n+1}^-\}.
\end{equation}
and indexed by
$$
E_{n+1}=(E_n\setminus \{\ve\})\cup \hbox{Succ}({\ve})\subset \eE^*.
$$

Define the sequence of jump times $\pi_n$ by 
$$
\forall n\in \NN^*:\quad \pi_n=\zeta (\Y_n)^{-1}V_n, 
$$
such that their conditional distributions satisfy,
\begin{equation}
\label{distpi}
\pi_n \, | \, \zeta (\Y_n) \sim \hbox{Exponential}(\zeta (\Y_n)).
\end{equation}
As $(V_n: n\in \NN^*)$ is a sequence of independent random variables,
$(\pi_n: n\in \NN^*)$ is conditionally independent given $\sigma(\Gamma_n: n\in \NN)$.
Define the sequence of times 
$$
\aS_0=0 \, \hbox{ and } \; \aS_n=\aS_{n-1}+\pi_{n}=\sum_{j=1}^{n} \pi_j \,
\hbox{ for } n\in \NN^*.
$$
We have
$\lim\limits_{n\to \infty}\aS_n=\infty\,$ $\PP-$a.s.
because $(N(t): t\ge 0)$ defined by 
\begin{equation}
\label{defNYY}
N(t)=\sup\{n\in \NN: \aS_n\le t), 
\end{equation}
is stochastically dominated by a birth chain 
$(M(t): t\ge 0)$, with $M(0)=1$
and linear birth rates $b_n=n\, \Lambda([W])$.

\medskip

Define the tessellation process $\wY \wedge W=((\wY\wedge W)_t: t\ge 0)$ 
taking values on $\T\wedge W$, by
\begin{equation}
\label{defYp}
(\wY\wedge W)_t=\Y_n  \, \hbox{ when } t\in
\left[\aS_{n}, \aS_{n+1}\right), n\in \NN,
\end{equation}
where $(\Y_n:n\in \NN)$ is the sequence of tessellations defined in the
algorithm. Hence, $\wY\wedge W$ is well-defined
for all times $t\ge 0\;$ $\,\PP-$a.s. We also have
$$
\forall n\in \NN: \quad (\wY\wedge W)_{\aS_n}=\Y_n,
$$
and so $(\aS_n: n\in \NN^*)$ is the sequence
of jump times of  $\wY\wedge W$.

\medskip

\begin{theorem}
\label{igual1}
$\wY\wedge W$ is the STIT process associated to $\Lambda_{[W]}$.
\end{theorem}

{\em Proof:}
The proof of the Markov property of $\wY\wedge W$ is analogous to
the one made in Proposition \ref{primeas}.

\medskip

The tessellation $\Y_{n+1}$ is formed from
$\Y_n$ by the division of the random cell $C^*_n$ of $\Y_n$, 
chosen with the help of $H_{n+1}$ and $U_{n+1}$, so
\begin{equation}
\label{reldep3}
C^*_n\in \sigma(\Y_n,H_{n+1},U_{n+1})
\end{equation}

\medskip

As said, to define $\Y_{n+1}$ we need to define the random cell
$C^*_n$ which is divided by hyperplane $H_{n+1}$. 
Hence, $\Y_{n+1}\in \sigma(\Y_n, U_{n+1},H_{n+1})$ 
and by recursion we get

\begin{equation}
\label{reldep5}
\Y_{n+1}\in \sigma(U_k,H_k: k\le n+1) .
\end{equation}
From (\ref{reldep2}), (\ref{reldep1}) and (\ref{reldep5}) 
and also by  recursion  we find,
\begin{equation}
\label{reldep5'}
H_{n+1}\in \sigma(\Y_n, G_{n+1})\subseteq 
\sigma(U_k : k\le n; G_k: k\le n+1).
\end{equation}
Since $\sigma(U_n, G_n: n\in \NN) \, \qperp \,  \sigma(V_n:n\in \NN)$, 
we obtain
\begin{equation}
\label{reldep6}
\sigma(\Y_n, U_n, G_n, H_n: n\in \NN) \, \qperp \,  (V_n:n\in \NN).
\end{equation}
Note also that (\ref{reldep5}) and (\ref{reldep5'}) imply that
\begin{equation}
\label{reldep4}
\Y_n\, \qperp \, G_{n+1} \hbox{ and } H_{n+1}\, \qperp \, U_{n+1} \, | \, \Y_n.
\end{equation}

\medskip

From $\pi_{n+1} \,|\, \Y_n \sim \,$ Exponential$((\zeta (\Y_n)))$
we get that the distribution of the holding time at $t$, noted $\wtau_t$,
satisfies,
\begin{equation}
\label{indht33}
\wtau_t \, | \, \Y_{N(t)} \sim \hbox{Exponential}(\zeta(\Y_n)).
\end{equation}
Now, from (\ref{reldep6}) and (\ref{distpi}) we get the conditional
independence relation  (\ref{indht1}), this is
$\pi_{n+1} \, \qperp \, \sigma(C^*_{n}, H_{n+1}) \; | \, \Y_n$.
So, 
\begin{equation}
\label{indht34}
\wtau_t \, \qperp \, \sigma(C^*_{N(t)}, H_{N(t)+1}) \; | \, \Y_{N(t)}.
\end{equation}

\medskip

\noindent (I) Let us prove that for ${\ve}\in E_n$, the
conditional probability given $\Y_n$, that $C^*_n=C(\ve)$ 
is the cell divided at time $\aS_{n+1}$ satisfies (\ref{newx2'}),
that is we must show that 
\begin{equation}
\label{lafor}
\PP(C^*_n=C(\ve) \, | \, \Y_n)=
\frac{\Lambda([C(\ve)])}{\zeta(\Y_n)}.
\end{equation}
Note that (\ref{forvol}) can be written 
$H_{n+1} \, | \, \Y_n \sim (\zeta(\Y_n))^{-1}
\Gamma_n$. Hence, by using this relation together with
(\ref{decx}) and (\ref{relcr2}) we get,
\begin{eqnarray}
\nonumber
\PP(C^*_n=C(\ve) \, | \, \Y_n)&=&
\zeta(\Y_n)^{-1} 
\int_{[W]} \PP\left(U_{n+1}\in [a^n_{\ve}, b^n_{\ve}) \, | 
\, H_{n+1}\right) d\Gamma_n(H_{n+1})\\
\nonumber
&=&\zeta(\Y_n)^{-1}
\int_{[W]} \xi_n(H_{n+1})^{-1}{\bf 1}_{\{H_{n+1}\in C(\ve)\}} d\Gamma_n(H_{n+1})\\
\nonumber
&=&\zeta(\Y_n)^{-1}\Lambda_{[C(\ve)]}([W])\\
\label{finI}
&=&\zeta(\Y_n)^{-1}\Lambda([C(\ve)]).
\end{eqnarray}
Hence (\ref{lafor}) follows.

\medskip

\noindent (II) We claim that
the distribution of the random hyperplane $H_{n+1}$, conditional
to $C^*_n=C(\ve)$ and $\Y_n$, is $\hLambda_{[C(\ve)]}$. 
To prove it we use that $H_{n+1}=\Xi_{n}\circ G_{n+1}$, 
that $\Xi_n$ only depends on $\Gamma_n$ and so on
$\Y_n$, and that 
$$
\Y_n\, \qperp \,  G_{n+1}, \; H_{n+1}\in \sigma(\Y_n, G_{n+1})
\hbox{ and } H_{n+1} \, \qperp \, U_{n+1} \, | \, \Y_n. 
$$
These relations allow to get for all $\K\in {\cal B}([W])$,
\begin{eqnarray}
\nonumber
&{}& \PP(H_{n+1}\in \K, C^*_n=C(\ve) \, | \, \Y_n)=
\int_\K \PP(dH_{n+1}, C^*_n=C(\ve) \, | \, \Y_n)\\
\nonumber
&{}& =
\int_\K \PP(dH_{n+1}, U_{n+1}\in [a^n_{\ve}, b^n_{\ve}) \, | \, \Y_n)
=\int_\K \xi_n(H_{n+1})^{-1}{\bf 1}_{\{H_{n+1}\in C(\ve)\}} \PP(dH_{n+1})\\
\label{conjH}
&{}& =\int_{C(\ve)\cap \K} \zeta(\Y_n)^{-1}d\Gamma_n(H)
=\zeta(\Y_n)^{-1}\Lambda_{[C(\ve)]}(\K)\\
\nonumber
&{}&=\zeta(\Y_n)^{-1}\Lambda([C(\ve)]\cap \K).
\end{eqnarray}
From (\ref{lafor}) we get  desired distribution:
\begin{equation}
\label{ultimcc}
\PP(H_{n+1}\in \K \, | \, C^*_n=C(\ve), \Y_n)=
(\Lambda([C(\ve)])^{-1}\Lambda([C(\ve)]\cap \K),
\end{equation}
Then, from (\ref{indht33}),
(\ref{ultimcc}) and (\ref{lafor}), we get that 
for all $C(\ve)\in \Y_{N(t)}$, $H\in [C(\ve)]$ and $s>0$ it is satisfied
$$
\PP(H(\ve)\in dH, C^*_t=C(\ve), \wtau_t\in ds  \, | \, \Y_{N(t)})=
\hLambda(dH)\, \Lambda([C(\ve)])\, e^{-\zeta(\Y_{N(t)})\, s}ds.
$$
We have shown (\ref{lawX0}) and so 
$\wY\wedge W$ given by (\ref{defYp}) is a STIT tessellation
associated to $\Lambda_{[W]}$.  \hfill $\Box$

\medskip

With respect to relation (\ref{relcr2}): we have first selected
an hyperplane with a probability measure $\hGamma_n$ (so depending
on the tessellation), and the 
equiprobability relation (\ref{relcr2}) for the cells 
which are intersected by the hyperplane, is nothing but 
an explicit computation of $(d \Lambda_{[C(\ve)]}/d\Gamma_n)(H)$.
The equiprobability relation appearing in \cite{st} page $9$,   
is in a different context and the hyperplane is chosen 
with probability measure $\Lambda$.

\bigskip

\noindent {\bf Example}.
As an illustrative example let us see what happens in the case $n=2$. 
The hyperplane $H_1$ divides $C$ in two cells $C(+)$ and $C(-)$,
and so $\Gamma_{1}=\Lambda_{[C(+)]}+\Lambda_{[C(-)]}$. 

\medskip

\noindent If $\{ H_2 \} \cap [C(+)]\cap [C(-)]=\emptyset$, then $C^*_1=C(e)$
when $H_2\in [C(e)]$.

\medskip

\noindent If $\{ H_2\} \cap  [C(+)]\cap [C(-)]\neq \emptyset$
then $\xi_1 (H)=\frac{1}{2}$. 
Following the order on $E_2$ the decision is: 
if $U_1\in [0,1/2)$ then $C^*_1=C(-)$ 
and if $U_1\in [1/2,1)$ then $C^*_1=C(+)$.
We note that
$[C(+)]\cap [C(-)]=[C(+)\cap C(-)]$,
and so $H_2\cap  [C(+)]\cap [C(-)]\neq \emptyset$ is equivalent to
$H_2\cap[C(+)\cap C(-)] \neq \emptyset$. 

\medskip

\subsection {Simulation of random hyperplanes}
\label{remsim}
Because it is not obvious from (\ref{reldep2}) how to generate $H_{n+1}$, we 
provide here a description of its density (\ref{relcr0}) which may be 
used in a simulation. Denote by $\pi_u C$ the orthogonal projection 
of $C$ onto the one-dimensional linear subspace (of $\RR^\ell$) 
spanned by $u\in {\SaS}_+^{\ell -1}$, and by $\lambda (\pi_u C)$ the 
length of this projection, which is also called the width or breadth of $C$ 
in direction $u$. If the image of $\Lambda$ on the parameter space is given 
by (\ref{prodmeas}), then the density  $\tilde \xi_n$ of the parametric 
representation of $H_{n+1}$ is given by
\begin{equation}
\label{densityhyper}
\tilde \xi_n (\alpha ,u)\, \lambda (d\alpha ) \, \theta (du) =
\frac{\displaystyle{\sum_{C\in \Y_n} 1_{\pi_u C}(\alpha )}}
{\displaystyle{\sum_{C\in \Y_n} \lambda (\pi_u C)}} \lambda (d\alpha ) \ 
\frac{\displaystyle{\sum_{C\in \Y_n}\lambda (\pi_u C)}}
{\displaystyle{\sum_{C\in \Y_n}\lambda (\pi_{u'} C)\theta (du')}} \theta (du).
\end{equation}
Of course, the sum ${\sum_{C\in \Y_n}\lambda (\pi_u C)}$ can be canceled out, 
but the given form shows better the decomposition of the joint density of the 
two parameters into a probability density w.r.t. $\theta$ for the direction 
$u\in {\SaS}_+^{\ell -1}$ and a conditional probability density, given $u$, 
for $\alpha \in \RR$, w.r.t. the Lebesgue measure. Note that 
$\sum_{C\in \Y_n} 1_{\pi_u C}(\alpha )=\xi_n(H(\alpha ,u))=
\#\{C\in \Y_n: H(\alpha ,u)\in [C]\}$ (see (\ref{relcr0})), 
and hence the conditional density of $\alpha$ is a step function.

\bigskip

\noindent {\bf Acknowledgments}.  The authors are indebted for the support 
of Program Basal CMM PFB03 from CONICYT-Chile and 
from DFG Germany (Project NA247/6-2).

\end{document}